\documentclass[12pt]{article}

\usepackage[cp1251]{inputenc}
\usepackage[T2A]{fontenc}
\usepackage{amsmath,amsfonts,amssymb,amsthm,cite}
\usepackage{geometry}
\theoremstyle{plain}
\newtheorem{theorem}{Theorem}
\newtheorem{lemma}{Lemma}

\newtheorem{definition}{Definition}
\theoremstyle{definition}

\newcommand{\nonprint}[1]{}

\begin{document}

\begin{flushleft}

\Large
\textbf{Vladimir Mikhailets,  Olena Atlasiuk, Tetiana Skorobohach}

\small

\vspace{+0.3cm}

mikhailets@math.cas.cz, hatlasiuk@gmail.com, tetianaskorobohach@gmail.com

\vspace{+0.5cm}

\Large

\textbf{On the solvability of Fredholm boundary-value problems\\ in fractional Sobolev spaces}
\end{flushleft}

\normalsize

\medskip

\begin{abstract}
Systems of linear ordinary differential equations with the most general inhomogeneous boundary conditions in fractional Sobolev spaces on a finite interval are studied. The Fredholm property of such problems in corresponding pairs of Banach spaces is proved, and their indices and dimensions of kernels and cokernels are found. Examples are given that show the constructive character of the obtained results.
\medskip
\end{abstract}

\section{Introduction and statement of the problems}\label{section1}

The investigation of solutions of systems of ordinary differential equations is an important part of numerous problems of contemporary analysis and its applications (see, e.g., monograph \cite{BochSAM2004} and the references therein). Unlike Cauchy problems, the solutions of such problems may not exist or may not be unique.

For inhomogeneous boundary-value problems on a finite interval of the form
\begin{equation*}
Ly:=y^{\prime}(t) + A(t)y(t)=f(t), \quad t\in(a,b),
\end{equation*}
\begin{equation*}
By=c,
\end{equation*}
where the matrix-valued function $A(\cdot)$ and the vector-valued function $f(\cdot)$ are summable on $[a, b]$, and the linear continuous operator
$$
  B\colon C \big([a,b];\mathbb{R}^{m}\big)  \rightarrow\mathbb{R}^{m},
$$
the questions of correct solvability and continuous dependence of solutions in a parameter in the space $C \big([a,b];\mathbb{R}^{m}\big)$ were studied in the papers of I.\,T.~Kiguradze \cite{Kigyradze1975, Kigyradze1987} and his followers \cite{KodliukMikhailets2013JMS, MPR2018, MikhailetsChekhanova}. Such problems cover all classical types of boundary conditions (two-point, multi-point, integral, mixed), but do not cover problems containing derivatives of an unknown function of integer or fractional orders in boundary conditions. Such boundary conditions are related to the function spaces in which the problem is studied.  Their analysis requires new research approaches and methods. In the case of Sobolev spaces of integer order, their analysis was carried out in \cite{GKM2015, KodlyukM2013, GKM2017, AtlMikh2019}, and in the case of H\"{o}lder spaces in the paper \cite{MMS2016}. At the same time, the analytical description of linear operators continuously acting from Sobolev space or $C^{(n)}$ into the space $\mathbb{C}^{m}$ was essentially used.

In this paper, the case of \textit{fractional} Sobolev spaces is investigated. For such spaces, there is no description of linear continuous operators acting from these spaces in $\mathbb{C}^{m}$, which significantly complicates the study of boundary-value problems.

Let's introduce some necessary notations for statement of the problem. Let the finite interval $(a,b)\subset\mathbb{R}$ and numerical parameters be given $$\{m, n, r\} \subset \mathbb{N}, \quad s\in (1,\infty)\setminus \mathbb{N}, \quad 1\leq p<{\infty}.$$

By $W_p^n:=W_p^n([a,b];\mathbb{C})$ we denote a complex Sobolev space and set  $W_p^{0}:=L_p$. We denote Sobolev spaces of vector-valued functions $(W_p^n)^{m}:=W_p^n([a,b];\mathbb{C}^{m})$ and matrix-valued functions$(W_p^n)^{m\times m}:=W_p^n([a,b];\mathbb{C}^{m\times m})$, respectively, whose elements belong to the function space $W_p^n$. By $\|\cdot\|_{n, p}$ we also denote the norms in these spaces. They are defined as the~sums of~the corresponding norms of all elements of a vector-valued or matrix-valued function in $W_p^n$. The~space of functions (scalar, vector-value, or matrix-value functions) in which the~norm is introduced is always clear from the context. For $m=1$, all these spaces coincide. As it is well known, the spaces $W_p^n$ are Banach and separable for $p<\infty$.

We denote by $W_p^s:=W_p^s([a,b];\mathbb{C})$, where $1\leq p<\infty$ and a non-integer $s>1$, Sobolev--Slobodetsky space of all complex-valued functions that belong to Sobolev space $W_p^{[s]}$ and satisfy the condition
$$
\|f\|_{s,p}:=\|f\|_{[s],p}+\left(\int\limits_a^b\int\limits_a^b\frac{\left|f^{[s]}(x)-f^{[s]}(y)\right|^p}{|x-y|^{1+\{s\}p}}dxdy\right)^{1/p}<+\infty.
$$
Here, $[s]$ is an integer, and $\{s\}$ is a fractional part of a number $s$. Here, we recall, $\|\cdot\|_{[s],p}$ is the norm in the Sobolev space  $W_p^{[s]}$. This equality defines the norm of $\|f\|_{s,p}$ in space $W_p^s$.

Consider on a finite interval $(a,b)$ a linear boundary-value problem for the system of $m$ differential equations of the first order
\begin{equation}\label{bound_pr_1}
(Ly)(t):=y^{\prime}(t) + A(t)y(t)=f(t), \quad t\in(a,b),
\end{equation}
\begin{equation}\label{bound_pr_2}
By=c,
\end{equation}
where the matrix-valued functions $A(\cdot)$ belong to the space $(W_{p}^{s-1})^{m\times m}$, the vector-valued function $f(\cdot)$ belongs to the space $(W_{p}^{s-1})^{m}$, the vector $c$ belongs to the space $\mathbb{C}^{r}$, and $B$ is a linear continuous operator
$$
  B\colon (W_{p}^{s})^{m} \rightarrow\mathbb{C}^{r}.
$$

The boundary condition \eqref{bound_pr_2} consists of $r$ scalar boundary conditions for system of $m$ differential equations of the first order. We represent vectors and vector-valued functions in the form of columns. In the case of $r>m$, the boundary-value problem \eqref{bound_pr_1}, \eqref{bound_pr_2} is \textit{overdetermined}, and for $r<m$ the problem is  \textit{underdetermined}. A solution to the~boundary-value problem \eqref{bound_pr_1}, \eqref{bound_pr_2} is understood as a vector-valued function $y(\cdot)\in (W_{p}^{s})^m$ satisfy\-ing equation \eqref{bound_pr_1} for $s>1+1/p$ everywhere, and for $s\leq 1+1/p$ almost everywhere on $(a,b)$, and equality \eqref{bound_pr_2} specifying $r$ scalar boundary conditions.

The solutions of equation \eqref{bound_pr_1} fill the space $(W_{p}^{s})^m$ if its right-hand side $f(\cdot)$ runs through the space $(W_{p}^{s-1})^m$. Hence, the~boundary condition \eqref{bound_pr_2} is the most general condition for this equation. It includes all known types of classical boundary conditions, namely, the Cauchy problem, two- and many-point problems, integral and mixed problems, and numerous nonclassical problems. The last class of problems may contain the derivatives integer or fractional order $\beta$ of the unknown vector-valued function, where $$0 \leq \beta < s-\frac{1}{p}.$$

The main aim of the present paper is to prove the Fredholm property for boundary-value problem \eqref{bound_pr_1}, \eqref{bound_pr_2} and to find its index. Moreover, we establish the dimensions of the kernel and cokernel
of the operator of inhomogeneous boundary-value problem due to similar properties of a special rectangular numerical
matrix. In the case of Sobolev spaces of integer order, similar results were obtained earlier in the paper \cite{AtlasiukMikhailets2019DNANU}.

\section{Main results}\label{section2}

We rewrite the inhomogeneous boundary-value problem \eqref{bound_pr_1}, \eqref{bound_pr_2} in the form of a linear operator equation
\[ (L,B)y=(f,c), \]
where $(L,B)$ is a linear operator in the pair of Banach spaces
\begin{equation}\label{th2-LB}
(L,B)\colon (W^{s}_p)^m\rightarrow (W^{s-1}_p)^m\times\mathbb{C}^r.
\end{equation}

Let $X$ and $Y$ be Banach spaces. A linear continuous operator $T\colon X \rightarrow Y$ is~called a Fredholm operator if its kernel $\ker T$ and cokernel $Y/T(X)$ are finite-dimensional. If~operator $T$ is Fredholm one, then its range $T(X)$ is closed in $Y$ and the index
$$
\mathrm{ind}\,T:=\dim\ker T-\dim\big(Y/T(X)\big)\in \mathbb{Z}
$$
is finite $\big($see, e.g., \cite[Lemma~19.1.1]{Hermander1985}$\big)$.

\begin{theorem}\label{th_fredh-bis}
The linear operator \eqref{th2-LB} is a bounded Fredholm operator with the index $m-r.$
\end{theorem}
We denote by $Y(\cdot) \in (W_p^{s})^{m \times m}$ the unique solution of a linear homogeneous matrix equation with Cauchy initial condition:
\begin{equation} \label{Koshi}
Y'(t)+A(t) Y(t)=O_m, \quad t \in (a,b), \quad Y(a)=I_m.
\end{equation}
Here, $O_m$ is the zero $(m \times m)$ -- matrix, and $I_m$ is the identity $(m \times m)$ -- matrix. The unique  solution of Cauchy problem \eqref{Koshi} belongs to the space $(W_p^{s})^{m \times m}$.

\begin{definition}\label{matrix_BY}
A bloc numerical matrix
\begin{equation}\label{matrutsa}
M(L,B) \in \mathbb{C}^{m \times r}
\end{equation}
is characteristic matrix for the boundary-value problem \eqref{bound_pr_1}, \eqref{bound_pr_2}, if its j-th column is the result of the action of the  operator $B$ on the j-th column of the matrix-valued function $Y(\cdot)$.
\end{definition}
Here, $m$ is the number of scalar differential equations of system \eqref{bound_pr_1}, and $r$ is the number of scalar boundary conditions.

\begin{theorem}\label{th dimker}
The dimensions of the kernel and cokernel of the operator \eqref{th2-LB} are equal to the~dimensions of the kernel and cokernel of the characteristic matrix, respectively, 
\begin{equation}\label{dimker}
\operatorname{dim} \operatorname{ker}(L,B)=\operatorname{dim} \operatorname{ker}\big(M(L,B)\big),
\end{equation}
\begin{equation}\label{dimcoker1}
\operatorname{dim} \operatorname{coker}(L,B)=\operatorname{dim} \operatorname{coker}\big(M(L,B)\big).
\end{equation}
\end{theorem}

Necessary and sufficient conditions for the invertibility of the operator $(L, B)$ follows from Theorem \ref{th dimker}, that is, the condition under which problem \eqref{bound_pr_1}, \eqref{bound_pr_2} possesses a unique solution and this solution depends continuously on the right-hand sides of the differential equation and boundary condition.

\begin{theorem}\label{th_invertible-bis}
The operator $(L,B)$ is invertible if and only if $r=m$ and the square matrix $M(L,B)$ is nondegenerate.
\end{theorem}

\section{Examples}\label{section3}

\textit{Example 1.} Let us consider a linear \textit{one-point} boundary-value problem for a differential equation
\begin{equation}\label{1.6.1t1}
    Ly(t):= y'(t)+Ay(t)=f(t),\quad
t \in(a,b),
\end{equation}
\begin{equation}\label{1.3t1}
By= \sum _{k=0}^{n-1} \alpha_{k} y^{(k)}(a)=c.
\end{equation}
\noindent Here, $A$ is the constant $(m \times m)$ -- matrix, the vector-valued  function $f(\cdot)$ belongs to the space~$(W_{p}^{s-1})^{m}$, matrices $\alpha_{k} \in \mathbb{C}^{r\times m}$, the vector $c \in \mathbb{C}^{r}$, linear continuous operators $$B\colon (W_{p}^{s})^{m} \rightarrow\mathbb{C}^{r}, \quad (L,B)\colon (W^{s}_p)^m\rightarrow (W^{s-1}_p)^m\times\mathbb{C}^r,$$ the vector-valued  function $y(\cdot)\in (W_{p}^{s})^m$, and $s> n+ \frac{1}{p}-1$.

We denote by $Y(\cdot)\in (W_p^s)^{m\times m}$ the unique solution of Cauchy matrix problem
\begin{equation*}\label{r31}
   Y'(t)+A Y(t)=O_{m},\quad t\in (a,b), \quad Y(a)=I_{m}.
  \end{equation*}

Then the matrix-valued function $Y(\cdot)$ and its $k$-th derivative will have the following form:
\begin{gather*}
Y(t)= \operatorname{exp}\big(-A(t-a)\big), \quad Y(a) = I_{m}; \\
Y^{(k)}(t)= (-A)^k \operatorname{exp}\big(-A(t-a)\big), \quad Y^{(k)}(a) = (-A)^k, \quad k \in \mathbb{N}.
\end{gather*}

Substituting these values into the equation \eqref{1.3t1}, we have
$$M(L,B)=\sum_{k=0}^{n-1}\alpha_{k}(-A)^k.$$

It follows from Theorem \ref{th_fredh-bis} that $ind(L, B)=ind(M(L, B))= m-r$.

Therefore, by Theorem \ref{th dimker}, we obtain
\begin{gather*}
\operatorname{dim} \operatorname{ker}(L,B)=\operatorname{dim} \operatorname{ker}\left(\sum_{k=0}^{n-1}\alpha_{k}(-A)^k\right)=m -
\operatorname{rank}\left(\sum_{k=0}^{n-1}\alpha_{k}(-A)^k\right),  \\
\operatorname{dim} \operatorname{coker}(L,B)=-m+r+\operatorname{dim} \operatorname{ker}\left(\sum_{k=0}^{n-1}\alpha_{k}(-A)^k\right)=r -
\operatorname{rank}\left(\sum_{k=0}^{n-1}\alpha_{k}(-A)^k\right).\label{dimcoker}
\end{gather*}

It follows from these formulas that the Fredholm numbers of the problem do not depend on the choice of the right end $b > a$.

\textit{Example 2.} Let us consider a \textit{two-point} boundary-value problem for the system of differential equations \eqref{1.6.1t1} with the coefficient $A(t) \equiv O_{m}$ and the boundary conditions at the points $\{t_0, t_1\} \subset [a , b]$ containing derivatives of integer and $/$ or \textit{fractional} orders (in the sense of Caputo, see, for example, \cite{Kilbas2006}). They are given by equality
\begin{equation*}\label{3.BY1}
By=\sum _{j} \alpha_{0j} y^{(\beta_{0j})}(t_{0})+\sum _{j} \alpha_{1j} y^{(\beta_{1j})}(t_{1}).
\end{equation*}
Here, both sums are finite, the numerical matrices $\alpha_{kj} \in \mathbb{C}^{r \times m}$, and the nonnegative numbers $\beta_{kj}$ are such that for all $k \in \{1, 2\}$
\begin{equation*}\label{3.BY1}
\beta_{k,0}=0, \quad \beta_{kj}<s- 1/p.
\end{equation*}
Then, as is easy to verify, the linear operator
$$B\colon (W_{p}^{s})^{m} \rightarrow\mathbb{C}^{r}$$
is continuous.

It follows from Theorem \ref{th_fredh-bis} that the index of the operator $(L, B)$ is equal to $m-r$. We find its Fredholm numbers. In this case, the matrix-valued function $Y(\cdot)=I_{m}$. Therefore, the characteristic matrix has the form
\begin{equation*}\label{3.BY1}
M(L, B)= \sum _{j} \alpha_{0j} I_{m}^{(\beta_{0j})}+\sum _{j} \alpha_{1j} I_{m}^{(\beta_{1j})}= \alpha_{0,0}+\alpha_{1,0},
\end{equation*}
because the derivatives $I_{m}^{(\beta_{kj})}=0$ if $\beta_{kj} >0$ \cite{Kilbas2006}. Therefore, according to Theorem \ref{th dimker},
\begin{gather*}
\operatorname{dim} \operatorname{ker}(L,B)=\operatorname{dim} \operatorname{ker}\left(\alpha_{0,0}+\alpha_{1,0}\right)=m -
\operatorname{rank}\left(\alpha_{0,0}+\alpha_{1,0}\right),  \\
\operatorname{dim} \operatorname{coker}(L,B)=- m+r+ \operatorname{dim} \operatorname{ker}\left(\alpha_{0,0}+\alpha_{1,0}\right)=r -
\operatorname{rank}\left(\alpha_{0,0}+\alpha_{1,0}\right).
\end{gather*}

It follows from these formulas that the Fredholm numbers of the problem do not depend on the choice of the interval $(a, b)$, the points $\{t_0, t_1\} \subset [a, b]$ and the matrices $\alpha_{0, j}$ , $\alpha_{1, j}$, with $j \geq 1$.

\section{Preliminary results}\label{section4}

To prove Theorems \ref{th_fredh-bis}, \ref{th dimker}, \ref{th_invertible-bis}, we will need two auxiliary statements.

Let us introduce the metric space of matrix-valued functions
$$
\mathcal{Y}_{p}^{s}:=\bigl\{Y(\cdot)\in (W_{p}^{s})^{m\times m} \colon \quad Y(a)=I_{m},\quad \det Y(t)\neq 0\bigr\},
$$
with metric
$$
d_{p}^{s}(Y,Z):=\|Y(\cdot)-Z(\cdot)\|_{s,p}.
$$

\begin{theorem}\label{th 2.1}
Nonlinear mapping $\gamma \colon A\mapsto Y$, where $A(\cdot)\in(W_{p}^{s-1})^{m\times m}$, and $Y(\cdot) \in (AC[a,b])^{m\times m}$ is the solution of Cauchy problem \eqref{Koshi}, is a homeomorphism of Banach space $(W_{p}^{s-1})^{m \times m}$ on a metric space $\mathcal{Y}_{p}^{s}$.
\end{theorem}

The proof of Theorem \ref{th 2.1} is given in the article \cite{MikSko2021}.

We put
\begin{equation}\label{3.BY}
[BY]:=\left( B \begin{pmatrix}
                                              y_{1,1}(\cdot) \\
                                              \vdots \\
                                              y_{m,1}(\cdot) \\
                                            \end{pmatrix}
\ldots
                                    B \begin{pmatrix}
                                              y_{1,m}(\cdot) \\
                                              \vdots \\
                                              y_{m,m}(\cdot) \\
                                            \end{pmatrix}\right)=M(L, B).
\end{equation}

\begin{lemma}\label{dija oper}
For an arbitrary matrix-valued function $Y(\cdot)\in (W^{s}_p)^{m\times m}$, a vector $q\in\mathbb{C}^m$, and linear continuous operator $B \colon (W^{s}_p)^{m\times m} \times\mathbb{C}^r$, the equality holds
\begin{equation*}\label{rivn_matruz}
B(Y(\cdot)q) = \left[BY\right]q,
\end{equation*}
where the matrix $\left[BY\right]$ is defined by the equality \eqref{3.BY}.
\end{lemma}

\textit{Proof}. Let $i\in \{1,2,\ldots,m\}$, $k\in \{1,2,\ldots,m\}$, $j\in \{1,2,\ldots, r\}$, the matrix-valued function $Y(\cdot)=\left(y_{ik}(\cdot)\right)$, and the column vector $q = \left(q_{i}\right)$. Let's denote by  $$\left(\alpha_{j}\right)= \left[BY\right] q \quad \mbox{and}  \quad \left(\beta_{j}\right)= B(Y(\cdot) q).$$ Let
$$
B(y_{k}(\cdot)) =: (c_{j}).
$$
When the operator $B$ acts on the matrix-valued function $Y(\cdot)$, we get the matrix
$$
\left[BY\right] = (c_{ji}).
$$
Then we will get
$$
(\alpha_{j}) = (c_{ji}) \cdot (q_{i})
= \left(\sum _i c_{ji} q_{j} \right).
$$
Therefore, an arbitrary element $\alpha_{i}$ has the form
$$
\alpha_{j} = \sum _j c_{ji} q_{j}.
$$
But the following equalities hold
$$
(\beta_{j}) = B\left((y_{ik}(\cdot))\cdot (q_{k}) \right) = B \left( \sum _k  y_{ik}(\cdot)  q_{k} \right)=$$
$$= \sum _k  \left(B y_{ik}(\cdot)\right) q_{k} =\sum _k \left( c_{jk} \right) q_{k}=\left(\sum _j c_{jk}q_{j} \right).
$$
It follows that $\alpha_{j} = \beta_{j}$.

The proof is complete.

\section{Proofs of Theorems}\label{section5}

\textit{Proof of Theorem \ref{th_fredh-bis}}. Let us first justify the continuity of the operator $(L,B)$.
Since the operator $B$ is linear and continuous by convention, it suffices to prove the continuity of the operator $L$, which is equivalent to its boundedness. Boundedness of the linear operator
$$
L\colon (W^{s}_p)^m\rightarrow (W^{s-1}_p)^m
$$
follows from the definition of norms in Sobolev spaces $W_p^{s-1}$
and the well-known fact that each of these spaces forms a Banach algebra.

Let us now prove that the operator $(L,B)$ is Fredholm one and find its index. Let us choose a fixed linear bounded operator $C_{r,m}\colon (W^{s}_p)^m\rightarrow\mathbb{C}^{r}$. The operator $(L,B)$ admits the representation
$$
(L,B)=(L,C_{r,m})+(0,B-C_{r,m}).
$$
Here, the operator
$$
(L,C_{r,m})\colon (W^{s}_p)^m\rightarrow (W^{s-1}_p)^m\times\mathbb{C}^r,
$$
and the second term is a finite-dimensional operator. From the Second Stability Theorem (see, for example, \cite[Section~3, \S~1]{Kato_book}) it follows that the operator $(L,B)$ is Fredholm one if the operator $(L,C_{r, m})$ is such and
$$\operatorname{ind}(L,B)=\operatorname{ind}(L,C_{r,m}).$$
Therefore, it suffices to prove that the operator $(L,C_{r,m})$ is Fredholm one and to find its index by properly choosing the operator $C_{r,m}$. For this, we will consider three cases.

1. Let $r=m$.
Let's put
\[C_{m,m}y:=(y_1(a),\dots , y_m(a)).\]

Let's find the null space and the range of values of this operator. Let $y(\cdot)$ belongs to $\operatorname{ker}(L,C_{r,m})$. Then $Ly=0$ and $C_{m,m}y=(y_1(a),\dots , y_m(a))=0$. It follows from the theorem on the uniqueness of the solution of Cauchy problem that $y(\cdot)= 0$. Therefore, $\operatorname{ker}(L,C_{m,m})=0$.

Let $h\in (W^{s-1}_p)^m\times\mathbb{C}^m$ and $c\in\mathbb{C}^m$ are chosen arbitrarily. It follows from Theorem~\ref{th 2.1} that there exists a vector-valued function $y(\cdot)\in (W^{s}_p)^m$ such that
$$Ly=h, \quad (y_1(a),\dots, y_m(a))=c.$$ Then $\operatorname{ran}(L,C_{r,m})=\left(W^{s-1}_p\right)^m\times\mathbb{C}^m$.

2. Let $r>m$. Let's put
\[C_{r,m}y:=(y_1(a),\dots, y_m(a), \underbrace{0,\dots ,0}_{r-m})\in\mathbb{C}^{r}.\]

Let's find the null space of the operator $(L,C_{r,m})$. Let $y(\cdot)$ belongs to $\operatorname{ker}(L,C_{r,m})$. Then $Ly=0$ and $(y_1(a),\dots , y_m(a))=0$. From the theorem on the uniqueness of the solution of Cauchy problem, we have $y(\cdot)= 0$.

We write the set of values of the operator $(L,C_{r,m})$ in the form of a direct sum of two subspaces
\[\operatorname{ran}(L,C_{r,m})=\operatorname{ran}(L,C_{m,m})\oplus (\underbrace{0,\dots ,0}_{r-m}).\]
But, as proved before, $\operatorname{ran}(L,C_{m,m})=(W^{s-1}_p)^m\times\mathbb{C}^m$.

Hence, $\operatorname{def} \operatorname{ran}(L,C_{r,m})=r-m$.

3. Let $r\textless m$. Let's put
\[C_{r,m}y:=(y_1(a),\dots , y_r(a))\in\mathbb{C}^{r}.\]

We will prove that
$$
\operatorname{dim} \operatorname{ker}(L,C_{r,m})=m-r,
$$
$$
\operatorname{def} \operatorname{ran}(L,C_{r,m})=0.
$$
Let $y(\cdot)$ belongs to $\operatorname{ker}(L,C_{r,m})$. Then $Ly=0$ and $(y_1(a),\dots, y_r(a))=0$. Let us consider the following $m-r$ Cauchy problems:
\begin{gather*}
Ly_k=0, \quad C_{m,m}y_k=e_k, \quad \mbox{where} \quad k\in \{r+1, r+2,\dots ,m\}, \\
e_k:=(0,\dots, 0, \underbrace{1}_{k}, 0, \dots ,0) \in {C}^{m}.
\end{gather*}
It follows from Theorem ~\ref{th 2.1} that the solutions of these problems are linearly independent and form a basis in the subspace  $\operatorname{ker}(L,C_{r,m})$.

The surjectivity of the operator $(L,C_{r,m})$ follows from the already proven surjectivity of the operator~$(L,C_{m,m})$.

Hence, in each of the three cases, the operator $(L,B)$ is a Fredholm operator with an index $m-r$.

The proof is complete.

\textit{Proof of Theorem \ref{th dimker}}. Let us show that the equality~\eqref{dimker} is valid. Let's introduce the following notations:
\begin{gather*}
\operatorname{dim} \operatorname{ker}(L,B)=n',\\
\operatorname{dim} \operatorname{ker}\left(M\big(L,B\big)\right)=n''.
\end{gather*}
We justify the fulfillment of equality
\begin{equation}\label{riv nn}
n'=n''.
\end{equation}

Let $\operatorname{dim} \operatorname{ker}(L,B)=n'$. Then there are $n'$ such linearly independent solutions of the homogeneous equation $(L,B)y=(0,0)$ that
$$y_k(\cdot)\in \operatorname{ker}(L,B) \Leftrightarrow \left(\exists \, q_k \in\mathbb{C}^{m}\colon y_k(t) = Y(t) q_k, \quad \left[BY\right]q_k=0\right),$$ according to Lemma \ref{dija oper}, where the vectors $q_k\neq0$, and $k\in \{1, \ldots, n'\}$.
This means that $r-n'$ columns of the matrix \eqref{matrutsa} are linearly dependent and $n'\leq n''$.

On the contrary, let $\operatorname{dim} \operatorname{ker}\left(M\big(L,B\big)\right)=n''$, then its $r-n''$ columns are linearly independent. This means that for some vectors $q_k\neq0$, $k\in \{1, \ldots, n'\}$,
\begin{equation*}
\left[BY\right]q_k=0.
\end{equation*}
Let's put $$y_k(\cdot):=Y(\cdot) q_k.$$ Then $y_k(\cdot)\neq0$, $Ly_k(\cdot)=0$ and
\[By_k(\cdot)=B(Y(\cdot)q_k)=\left[BY\right]q_k=0,\]
based on Lemma \ref{dija oper}. Therefore, $y_k(\cdot)\in \operatorname{ker}(L,B)$, then $n'\geq n''$. Hence, the equality~\eqref{dimker} holds.

According to the definition, the characteristic matrix $M(L,B)$ belongs to the space $\mathbb{C}^{m\times r}$. As it is well known, the dimension of the kernel of the matrix is the difference between the number of its rows and its rank. And the dimension of the cokernel of the matrix is the difference between the number of columns and the rank. Then for the matrix $M(L,B)$, we have equality
\begin{equation}\label{coker M}
\operatorname{dim} \operatorname{coker}\left(M\big(L,B\big)\right)=r-m+\operatorname{dim} \operatorname{ker}\left(M\big(L,B\big)\right).
\end{equation}
From the formula for finding the index for the operator $(L,B)$
\begin{equation*}\label{coker LB1}
\mathrm{ind}\,(L,B):=\dim\ker(L,B)-\operatorname{dim} \operatorname{coker}(L,B),
\end{equation*}
we have
\begin{equation}\label{coker LB}
\operatorname{dim} \operatorname{coker}(L,B)=r-m+\dim\ker(L,B).
\end{equation}
The equalities \eqref{riv nn}, \eqref{coker M}, and \eqref{coker LB} imply the equality \eqref{dimcoker1}.

The proof is complete.

\section{Acknowledgements}\label{section6}

The research of the authors V. Mikhailets and O. Atlasiuk was supported by the  Czech Academy of Sciences, grant RVO:67985840.

\vspace{+0.2cm}

The research of the author O. Atlasiuk was supported by research work of young scientists of the National Academy of Sciences of Ukraine, 0121U111949, and research project of joint teams of scientists of Taras Shevchenko KNU and the National Academy of Sciences of Ukraine, 3M 2022.

\bigskip

\end{document}